\documentclass[letterpaper,12pt]{amsart}
\usepackage{amssymb}
\usepackage{xy}
\usepackage{epsfig}
\xyoption{all}

\newtheorem{theo}{Theorem}[section]
\newtheorem{lemma}[theo]{Lemma}
\newtheorem{coro}[theo]{Corollary}
\newtheorem{prop}[theo]{Proposition}
\newtheorem{conj}[theo]{Conjecture}
\newtheorem{defi}[theo]{Definition}
\theoremstyle{definition}
\newtheorem{rem}[theo]{Remark}

\def\GL{\operatorname{GL}\nolimits}
\def\ker{\operatorname{ker}\nolimits}
\def\Aut{\operatorname{Aut}\nolimits}
\def\Id{\operatorname{Id}\nolimits}

\def\S{{\mathfrak{S}}}

\def\bB{{\mathbf{B}}}
\def\bw{{\mathbf{w}}}

\def\M{{\mathcal{M}}}
\def\P{{\mathcal{P}}}

\def\CM{{\mathbf{C}}}
\def\NM{{\mathbf{N}}}
\def\RM{{\mathbf{R}}}

\def\ZM{{\mathbf{Z}}}         

\title{Springer theory in braid groups and the Birman-Ko-Lee monoid}
\author{David Bessis, Fran\c cois Digne, Jean Michel}
\address{Yale University, Department of Mathematics, PO Box 208 283,
New Haven CT 06520-8283, USA}
\email{david.bessis@yale.edu}
\address{LAMFA, Universit\'e de Picardie-Jules Verne
33, Rue Saint-Leu 80039 Amiens Cedex France}
\email{digne@u-picardie.fr}
\address{LAMFA, Universit\'e de Picardie-Jules Verne
33, Rue Saint-Leu 80039 Amiens Cedex France}
\address{ Institut de Math\'ematiques, Universit\'e Paris VII,
  175, rue du Chevaleret 75013 Paris France}
\email{jmichel@math.jussieu.fr}

\begin{document}
\maketitle    

\markboth{\sc D. Bessis, F. Digne, J. Michel}{\sc Springer theory in braid
groups
and the Birman-Ko-Lee monoid}
\pagestyle{myheadings}

{\flushleft {\bf Abstract} \small
We state a conjecture about centralizers of certain roots
of central elements in braid groups, and check it for Artin braid groups
and some other cases. Our proof makes use of results by Birman-Ko-Lee.
We give a new intrinsic account of these results.
}            

\bigskip

{\bf Notations}

If $G$ is a group acting on a set $X$, we denote by $X^G$ the
subset of $X$ of elements fixed by all elements of $G$. If
$(X,x)$ is a pointed topological space, we denote by
$\Omega(X,x)$ the corresponding loop space, by $\sim$ the homotopy relation on
$\Omega(X,x)$ and by $\pi_1(X,x)$ the fundamental group. For all $n\in
\NM$,
we denote by $\mu_n$ the set of $n$-th roots of unity in $\CM$.

\smallskip
{\bf Introduction}

Let $V$ be a finite dimensional complex vector space and let $W\subset\GL(V)$
be an irreducible finite group generated by complex reflections (that is,
elements $s\in\GL(V)$ such that $\ker(s-\Id)$ is a hyperplane). We denote by
$\M$ the complement in $V$ of the union of all reflecting hyperplanes. The space
$\M/W$ is called the {\em complement of the discriminant of $W$}; its
fundamental group is, by definition, the braid group $\bB$
associated to $W$.

Let $\zeta\in\CM$ be a root of unity
of order $d$. Suppose $d$ is a regular number for $W$, as
defined in \cite{springer}.
Thus there exists $w\in W$ such that
$$\ker(w-\zeta\Id) \cap \M \neq \emptyset.$$
Let us denote by $\M(w)$ this space.
According to \cite{springer}, $C_W(w)$ is a complex reflection
group on $\ker(w-\zeta\Id)$, with hyperplanes the traces of those of $W$ on
this space.
As noticed in \cite{brmi}, if the basepoint $x$ chosen for $\M/W$ is in
the image of $\M(w)$, there is a natural way of lifting
$w$ to an element $\bw$ of $\bB$, and the braid group $\bB_{C_W(w)}$ is
$\pi_1(\M(w)/C_W(w))$.
Following the implicit ideas behind question 3.5 in \cite{brmi}, 
let us state the following
conjecture, which claims that Springer theory of regular elements can
be lifted to braid groups:

\begin{conj}
\label{laconjecture}
The natural morphism $\bB_{C_W(w)} \rightarrow \bB$ induces an isomorphism
between $\bB_{C_W(w)}$ and the centralizer $C_\bB(\bw)$.
\end{conj}

The elements $\bw$ as above are
the $d$-th roots of a generator of the center of the pure braid
group (which is the kernel of the map $\bB\to W$, and has a cyclic center),

Conjecture \ref{laconjecture} can be reformulated in a more intrinsic way:

Let $(X,x)$ be a pointed topological space and $G$ be a group 
acting on $(X,x)$ (by morphisms in the category of pointed
spaces, so $x\in X^G$).
The action of $G$ can be naturally extended to the loop space $\Omega(X,x)$.
As clearly $\forall g \in G, \forall \gamma,\gamma'\in\Omega(X,x),
\gamma\sim\gamma' \Leftrightarrow g(\gamma)\sim g(\gamma')$, this
 induces a morphism $G \rightarrow \Aut(\pi_1(X,x))$.
For a given $G$, the construction of $\pi_1(X,x)^G$ from $(X,x)$ is actually 
functorial, from the category of pointed topological spaces with
$G$-action to the category of groups. 
Thus the natural injection $(X^G,x)\subseteq(X,x)$ induces
a natural morphism from $\pi_1(X^G,x)=\pi_1(X^G,x)^G$ to $\pi_1(X,x)^G$.

If $W$ is as above a finite irreducible
complex reflection group in $\GL(V)$, we denote by $X$ the corresponding
space $\M/W$. For all integers $m$,
the inclusion $\mu_m\subset\CM^\times$ and the identification
of $\CM^\times$ with the center of $\GL(V)$ define a
natural quotient action of $\mu_m$ on $X$.
As noticed in \cite{bessis} 1.2, when $d$ is regular and $w$ and
$\bw$ are as defined above,
$\M(w)/W(w)$ is homeomorphic to its image $X^{\mu_d}$ in $X$.
Let $x\in X^{\mu_d}$.
It is an easy calculation to check that
 the subgroup of $\Aut(\pi_1(X,x))$
generated by the conjugation by $\bw$ coincides
with the one arising from the action of $\mu_d$ on $X$.

Thus the above conjecture is equivalent to statement
that the natural morphism
$$\pi_1(X^{\mu_d},x)\rightarrow
\pi_1(X,x)^{\mu_d}$$ is an isomorphism.

The fact that the morphism is an isomorphism doesn't depend
on the choice of $x\in X^{\mu_d}$.

Our main theorem checks the conjecture for some specific groups
(we use the standard notations from \cite{shto}):

\begin{theo}
\label{letheoreme}
Let $W$ be an irreducible complex reflection group. 
Suppose $W$ is of one of the following types: $\S_n$;
$G(p,1,r)$ with $p > 1$; $G_{4}$; $G_5$; $G_8$; $G_{10}$;
$G_{16}$; $G_{18}$; $G_{25}$; $G_{26}$;
$G_{32}$. Let $X$ be complement of the discriminant of $W$. 
Let $d$ be a regular number for $W$, and let $x\in X^{\mu_d}$.
The natural morphism $$\pi_1(X^{\mu_d},x)
\rightarrow \pi_1(X,x)^{\mu_d}$$
 is an isomorphism.
\end{theo}

\section{The local monoid}
\label{localmonoid}

We will deduce our theorem from the particular case where $W$ is the
symmetric group $\S_n$, and $\bB$ the Artin braid group on $n$ strings.
In \cite{BKL}, Birman, Ko and Lee describe a remarkable monoid for
this group. The properties of their monoid will
be crucial in our proof. However, contrary to what is done in
\cite{BKL}, where the Artin braid group is initially
given by the Artin presentation, we prefer to use its more intrinsic
definition as a fundamental group. Of course, both viewpoints give
(non canonically) isomorphic groups; however, we believe our reformulation
is more natural than the original description. 

Let $X_n$ be the space of subsets of $\CM$ of cardinal $n$, with
its natural topology. In the setting of the introduction,
when taking the natural irreducible reflection representation
of the symmetric group $W=\S_n$ on $\CM^{n-1}$,
the complement of the discriminant $\M/W$ is
homotopy equivalent (in a way
compatible with the action of $\CM^\times$) to $X_n$.
We choose the usual direct (= anti-clockwise) orientation on $\CM$. 

Let us choose a basepoint $x\in X$, and let $B_x = \pi_1(X_n,x)$. 
We define in this section a monoid $M_x$,
which is a set of group generators for $B_x$. The structure of $M_x$
depends on the choice of $x$: when $x$ is taken to be the ``usual'' basepoint
$\{1,\ldots,n\}$, the monoid $M_x$ will be isomorphic to the usual Artin monoid;
choosing $\mu_n$ (the set of $n$-roots of unity) will yield what we call the
Birman-Ko-Lee monoid, which is isomorphic to the one described in \cite{BKL}.

If $\gamma\in\Omega(X_n,x)$ and $z\in x$,
we denote by $\gamma_z$ the ``string'' of $\gamma$ with origin
$z$. It is a path
$[0,1] \rightarrow \CM$, with $\gamma_z(0)=z$ and $\gamma_z(1)\in x$.
The element $\gamma$ is uniquely determined by $(\gamma_z)_{z\in x}$.
Conversely, a set of $n$ such strings which do not intersect define an element
of $\Omega(X_n,x)$. 

\begin{defi}
A pair $\{z,z'\} \subset x$ is said to be {\em non obstructing} if and only
if the closed line segment $[z,z']$ intersects $x$ only at $z$ and $z'$.
\end{defi}

We denote by $S_x$ the set of non obstructing pairs of elements of $x$.

{\bf \flushleft Examples.} The two crucial examples are
$x = \{ 1, \ldots, n \}$ and $x=\mu_n$. The corresponding $S_x$
have respectively cardinal $n-1$ and $n(n-1)/2$.

\def\de(#1,#2){\delta_{\{#1,#2\}}}
\def\dst{\de(z,z')}   

For each $\{z,z'\} \in S_x$, we 
denote by $\dst\in B_x$ the generator of the monodromy naturally
associated to $[z,z']$, as in the appendix of \cite{bessis}. A representative
of this element can be defined for instance by the set of strings
$\gamma_{z''}(t)=z''$
 if $z''\notin\{z,z'\}$, $\gamma_z(t)=f(z,z',\epsilon)(t)$,
$\gamma_{z'}(t)=f(z',z,\epsilon)(t)$ for $\epsilon\in\RM^*_+$
small enough (depending only on $x$) where
$$\forall t, f(z,z',\epsilon)(t)= \frac{z+z'}{2} + \frac{z-z'}{2}\cos(\pi t)
+i \epsilon \frac{z-z'}{2} \sin(\pi t)$$
(the corresponding arc is a half-ellipse with great axis $[z,z']$ and
small axis of length
$\epsilon |z-z'|$. The exact choice of $f(z,z',\epsilon)$ is
not important. One could for example replace the
half-ellipse by a half-rhomb, or any other variation. However, for later
use where defining a loop up to homotopy will not be sufficient, it is
convenient for us to define this explicit element of $\Omega(X_n,x)$).   

\begin{defi}
The submonoid $M_x$ of $B_x$ generated by  
$$\{ \dst |  \{ z,z'\} \in S_x \}$$ is called the {\em local monoid} at $x$.
\end{defi}

Let $l$ be the natural length function on $B_x$ (the map $B_x\to\ZM$
induced by the discriminant application $\M\to\CM^\times$). The
generators of $S_x$, being generators of the monodromy, have length 1,
so the monoid $M_x$ is $\NM$-graded (only the trivial element being of
length 0).

We denote by $\prec$ the left divisibility relation in $M_x$, i.e.
$$\forall m,m'\in M_x, m\prec m' \Leftrightarrow \exists m'' \in M_x,
m m''=m'.$$ It results from the $\NM$-grading of $M_x$ that the relation
$\prec$ is a partial order on $M_x$.

Clearly, the application $ \{ z,z'\} \mapsto \dst $ is injective, so we
may identify $S_x$ with its image in $M_x$.
We will later study extensively the structure of $M_{\mu_n}$.
We complete this section with some first properties which are valid for
all $x$.

\begin{prop}
\label{Sgenerates}
The set $S_x$ is a set of group generators for
$B_x$ (i.e., by taking $S_x \cup \{ s^{-1} | s\in S_x \}$, one has a set
of monoid generators for $B_x$).
\end{prop}

\begin{proof}
Distinguish one point $z\in x$. Draw the segments $\{ [z,z'] | z' \in x \}$.
Some of them may be obstructing, but by splitting these into 
smaller ones, one gets a planar graph
connecting all points in $x$ and whose edges are non obstructing.
The result then follows from the main theorem
in \cite{sergiescu} and its reformulation in the appendix of \cite{bessis}.
\end{proof}

{\bf \flushleft Notations and conventions.} 
We write $\lambda \vdash x$ to say that $\lambda$ is
a partition of $x$, in the usual set theoretical sense. 
If $y\subset x$ and $\lambda \vdash y$,
we will also use $\lambda$ to denote the partition of $x$ obtained by 
completing $\lambda$ with parts of cardinal $1$. In other words, we will
sometimes, for convenience, omit to write the cardinal $1$ parts of a
 partition.
By {\em convex polygon}, we mean either a point (if the number of vertices
is $1$), a segment (if the number of vertices is
$2$), or (if the number of vertices if $3$ or more ) a non-degenerate
convex polygon, i.e. such that three vertices never lie on a same line.
When $y\subset \CM$, we write $\overline{y}$ for the convex hull of $y$.

We now extend the notion of {\em non obstructivity} to partitions of
$x$:

\begin{defi}
\begin{itemize}
\item A finite non-empty 
subset $y\subset \CM$ is said to be {\em convex} if and only if it is 
the set of vertices of a convex polygon in $\CM$.
\item A partition $\lambda$ of $x$ is said to be {\em non obstructing} if and
only if it satisfies the following two properties:
\begin{itemize}
\item Every part $\nu$ in $\lambda$ is convex.
\item If two parts $\nu_1$ and $\nu_2$ are distinct, their convex hulls
$\overline{\nu_1}$ and $\overline{\nu_2}$ do not
 intersect.
\end{itemize}
\item The set of non obstructing partitions of $x$ is denoted by $\P_x$.
We use the notation $\lambda \models x$ to express that $\lambda\in\P_x$.
\end{itemize}
\end{defi}

We had already a notion of non obstructing pairs. 
The partition $\{ \{  z,z'\} \} $ is non obstructing if and only if 
$\{ z,z'\}$ is a non obstructing pair.

Let $\nu$ be a part of a non obstructing partition of $x$.
Choose $z_1$ an arbitrary element of $\nu$.
For $i\in \{1,\ldots,k\}$, denote by $z_i$ the $i$-th element of $\nu$
for the {\bf clockwise} order on $\nu$ starting at $z_1$.
In the next lemma, when $k=1$, the product is the empty product, thus the
trivial element  in $M_x$.

\begin{lemma}\label{tourner}
The element
$\de(z_1,z_2)\de(z_2,z_3)\ldots \de(z_{k-2},z_{k-1})\de(z_{k-1},z_{k})\in M_x$
does not depend on the choice of $z_1$ in $\nu$.
\end{lemma}

\begin{proof}
This is a consequence of one of the Sergiescu relations (see
\cite{sergiescu}, 1.1(ii)
or \cite{bessis}, th\'eor\`eme A.6) namely that
$$\de(z_1,z_2)\de(z_2,z_3)\ldots \de(z_{k-2},z_{k-1})\de(z_{k-1},z_{k})=
\de(z_k,z_1)\de(z_1,z_2)\de(z_2,z_3)\ldots \de(z_{k-2},z_{k-1}).$$
The proof of this relation is by induction on $k$. The case $k=3$
is checked by a direct computation, and for other $k$ we have
$$\displaylines{
\de(z_1,z_2)\de(z_2,z_3)\ldots \de(z_{k-2},z_{k-1})\de(z_{k-1},z_k)
\hfill\cr\hfill=
\de(z_2,z_3)\ldots \de(z_{k-2},z_{k-1})\de(z_{k-1},z_1)\de(z_{k-1},z_k)
\cr\hfill=
\de(z_2,z_3)\ldots \de(z_{k-2},z_{k-1})\de(z_{k-1},z_k)\de(z_k,z_1)
\cr}$$
where the first equality is by induction and the second by
the case $k=3$.
\end{proof}

We denote $\delta_\nu$ the element of \ref{tourner} (when $\{ z,z'\}$ is
a non obstructing pair, $\delta_{\{ \{ z,z'\} \} }$ coincides with the element 
$\dst$ defined earlier).

\begin{lemma}
Let $\nu$ and $\nu'$ be two finite non-empty subsets of $\CM$. Suppose
$\{\nu,\nu'\}$ can be completed to a non obstructing partition
 of $x$. Then we have $\delta_\nu \delta_{\nu'} = \delta_{\nu'} \delta_\nu$.
\end{lemma}

\begin{proof}
The convex hulls of $\nu$ and $\nu'$ cannot intersect, thus the 
generators corresponding to their edges commute pairwise.
\end{proof}

Let $\lambda \models x$. The above lemma makes it natural (and non-ambiguous)
to define
$$\delta_{\lambda} = \prod_{\nu\in\lambda} \delta_\nu.$$

\begin{defi} Given two partitions $\lambda$ and $\lambda'$ of $x$, we say that 
$\lambda$ is {\em finer} than $\lambda'$ (or equivalently that 
$\lambda'$ is {\em coarser} than $\lambda$), and we write
$\lambda \prec \lambda '$, if and only if 
$\forall \nu \in \lambda, \exists \nu' \in \lambda', \nu \subset \nu'$.
\end{defi}

Clearly $\prec$ is a partial order on the set of partitions of $x$.
We intentionally use the same symbol $\prec$, as for the left divisibility order
in $M_x$. The next proposition justifies this notation. 
>From the definition of a non obstructing partition, it is clear that
if $\lambda$ is finer than $\lambda'$, and $\lambda' \models x$, then
$\lambda\models x$. 

Let $P_x = \{ \delta_\lambda | \lambda \in \P_x \}$.
The main result in this section is the following proposition.

\begin{prop}
\label{poset}
\begin{itemize}
\item[(i)]The left divisibility order in $P_x$ coincides with the restriction of
$\prec$ from $M_x$ to $P_x$; in other words,
$$\forall p,p' \in P_x, (\exists p'' \in P_x, p p''= p') \Leftrightarrow
(\exists m'' \in M_x, p m''= p').$$
We denote this partial order by $\prec$.

\item[(ii)] The map 
\begin{eqnarray*}
D: \P_x & \longrightarrow & P_x \\
\lambda & \longmapsto & \delta_\lambda
\end{eqnarray*}
is a poset isomorphism from $(\P_x,\prec)$ to
$(P_x,\prec)$. 
\end{itemize}
\end{prop}

Before proving the proposition, we need some definitions and lemmas.

\begin{defi}
\label{coupe} Let $\nu$ be a convex subset of $\CM$, and
let $\nu'$ be a subset of $\nu$.
 Let $z_1,\ldots,z_k$ be
a clockwise numbering of the elements of $\nu$, such that
$z_{i_1},\ldots,z_{i_{k'}}$ is a clockwise numbering of $\nu'$ with
$1=i_1<i_2<\ldots<i_{k'}\le k$.
Then we denote by $\nu'\backslash\nu$ ($\nu$ ``cut at'' $\nu'$) the
partition with parts the sets $\{z_{i_j},z_{i_j+1},\ldots,z_{i_{j+1}-1}\}$
(where for $j=k'$ we make the
convention that $i_{k'+1}=k+1$).
\end{defi}

{\bf \flushleft Example.} In the picture below, the grey points are
the points of $\nu'$, the other points of $\nu$ are black, and
the parts of $\nu'\backslash\nu$ are enclosed by dashed curves.
$$\epsfig{file=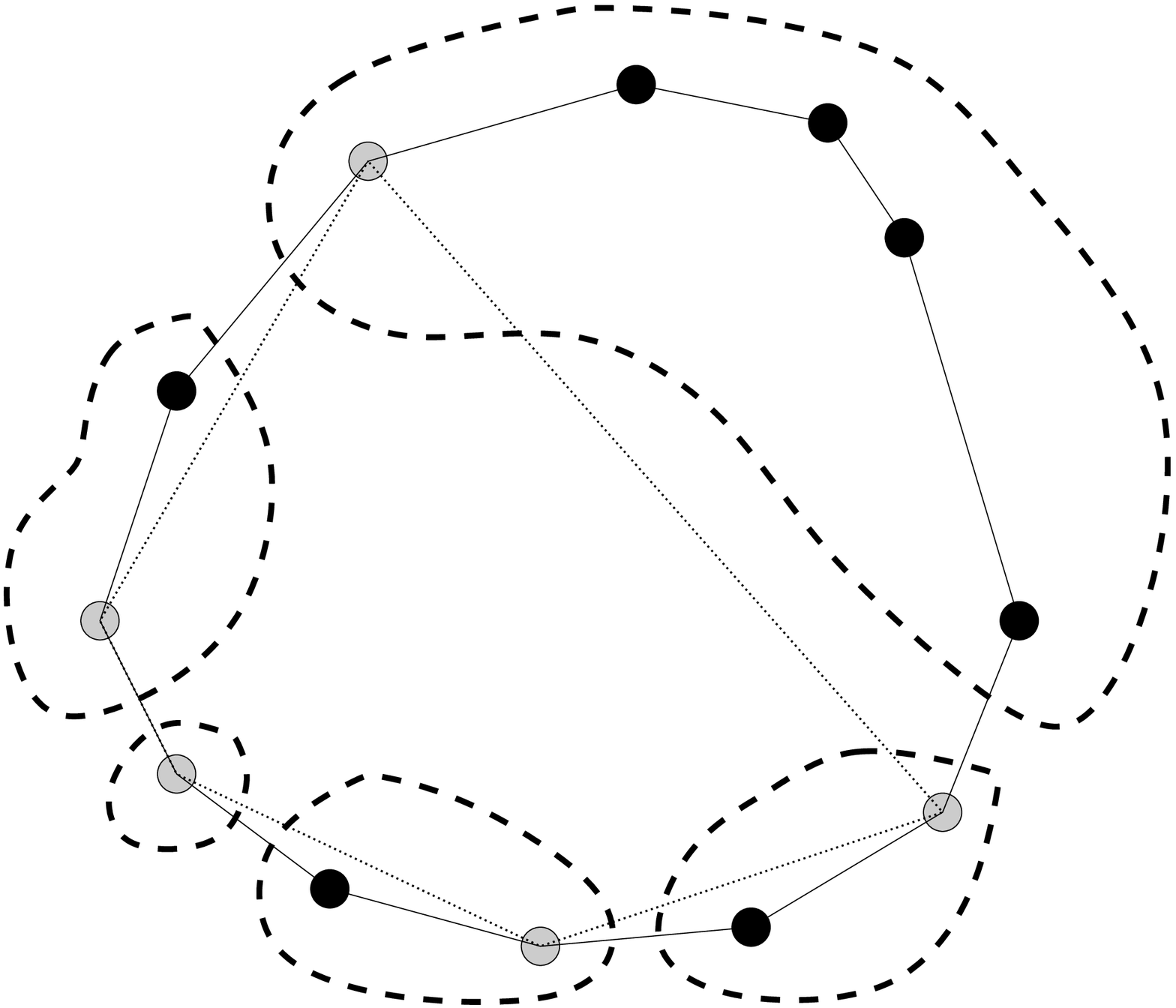,height=4cm}$$  

Note that $\nu'\backslash\nu$ is well defined (does not depend on the
chosen numbering) and is a non obstructing partition of $\nu$. We will
need the following alternative description of $\nu'\backslash\nu$:
its parts of cardinal $\geq 2$
are the intersections of $\nu$ with each connected component of
the complement of $\overline{\nu'}$ in $\overline\nu$, to which has been
added the element of $\nu'$ just before the connected component (in
clockwise order).

Denote by $\phi$ the natural epimorphism from $B_x$ to $\S_x$.
It maps $\dst$ to the transposition $(z,z')$. Thus,
when $\lambda \models x$, it is clear
by construction that $\lambda$ is the orbit decomposition for the action
in $x$ of $\phi (\delta_\lambda)$. Thus if $\lambda \neq \lambda'$,
then $\phi (\delta_\lambda)$ and $\phi (\delta_{\lambda'})$ have different cycle
decompositions, and are different. We have proved the:

\begin{lemma}
\label{phiinjective}
The restriction of $\phi$ to $P_x$ is injective.
\end{lemma}

\begin{lemma}\label{diviseurs}
\begin{itemize}
\item[(i)]  Let $\nu$ be a convex subset of $\CM$, and
let $\nu'$ be a subset of $\nu$.
We have
$\delta_\nu=\delta_{\nu'}\delta_{\nu'\backslash\nu}$.
\item[(ii)]
Suppose $\lambda,\lambda'\models x$,
 $\lambda'$ has only one part of cardinal $\geq 2$, and
$\lambda$ is finer that $\lambda'$.
Then there exists a unique non obstructing partition of $x$, which
we denote by $\lambda \backslash \lambda'$ such that
$\delta_\lambda \delta_{\lambda \backslash \lambda'}=\delta_{\lambda'}$.
\item[(iii)]
Suppose $\lambda,\lambda'\models x$ and $\lambda$ is finer that $\lambda'$.
Then there exists a unique non obstructing partition of $x$, which
we denote by $\lambda\backslash\lambda'$ such that
$\delta_\lambda \delta_{\lambda \backslash \lambda'}=\delta_{\lambda'}$.
\end{itemize}
\end{lemma}
\begin{proof}
We prove (i) by induction on the cardinality of $\nu'$.
With the notations of \ref{coupe},
let $\nu_1=\{z_1,\ldots,z_{i_2-1}\}$, $\nu_2=\nu-\nu_1$ so that
$\delta_\nu=\delta_{\nu_1}\de(z_{i_2-1},z_{i_2})\delta_{\nu_2}$, and
let $\nu'_2=\nu'-\{z_{1}\}$, so that
$\delta_{\nu'}=\de(z_1,z_{i_2})\delta_{\nu'_2}$.
Then by \ref{tourner} applied to $\nu_1\cup\{z_{i_2}\}$ we have
$\delta_\nu=\de(z_1,z_{i_2})\delta_{\nu_1}\delta_{\nu_2}$. 
By induction hypothesis, we have
 $\delta_{ \nu_2 } =\delta_{\nu_2'}\delta_{\nu_2'\backslash\nu_2}$.
As $\nu'\backslash\nu = \{ \nu_1 \} \cup \nu_2'\backslash\nu_2$,
and as $\delta_{\nu'_2}$ commutes to $\delta_{\nu_1}$ (since
$\overline{\nu'_2}\cap\overline{\nu_1}=\emptyset$),
we get the
result by induction.

(ii): denote by $\nu$ the only non trivial part of $\lambda'$, and
by $\lambda_1,\ldots,\lambda_l$ the non trivial parts of $\lambda$.
To prove the result, the essential step
is to notice that since
for $i>1$ we have
$\overline{\lambda_i}\cap\overline{\lambda_1}=\emptyset$, each
$\lambda_i$ lies inside a single connected component of the complement of
$\overline{\lambda_1}$ in $\overline\nu$, thus each is included in a part of
$\lambda_1\backslash\nu$. Thus (ii) follows by induction on the number
of parts of $\lambda$ from (i). The uniqueness comes from the fact that the
identity
$\delta_\lambda \delta_{\lambda \backslash \lambda'}=\delta_{\lambda'}$
is valid in $M_x\subset B_x$, and $B_x$ is a group.
 
Finally, (iii) is easily obtained by applying (ii) to all the non trivial parts
of $\lambda'$.
\end{proof}                                       

For $\sigma \in \S_x$, we denote by $|\sigma |$ the minimum
number of transpositions in a decomposition of $\sigma$ as a product of
transpositions; such a decomposition is called {\em reduced} (note that
we allow all transpositions). The map $\sigma \mapsto |\sigma|$ is not a
morphism, but we have the relation $|\sigma_1\sigma_2 | \leq
|\sigma_1| + |\sigma_2|$

\begin{lemma}
\label{length}
\begin{itemize}
\item[(i)] Let $\lambda \vdash x$ be the orbit
decomposition of $\sigma \in \S_x$. We have
$$|\sigma | = \sum_{\nu\in\lambda} (|\nu | -1) = |x| - |\lambda |$$
\item[(ii)] Let $\sigma\in\S_x$. Let $t=(z, z')$. If $z$ and $z'$ are
in the same orbit of $\sigma$, then $|\sigma t| = |  \sigma | - 1$.
If $z$ and $z'$ are
in different orbits of $\sigma$, then $|\sigma t| = | \sigma | + 1$.
\item[(iii)] Let $\sigma\in\S_x$. Let $t_1t_2\ldots t_{|\sigma|}$ be a reduced
decomposition of $\sigma$. Then if $t_i=(z, z')$ is one of the transpositions,
$z$ and $z'$ belong to the same orbit of $\sigma$. Let $\lambda \vdash x$
be the orbit decomposition of $\sigma$ and, for $i=1,\ldots,|\sigma |$, 
$\lambda_i \vdash x$
be the orbit decomposition of $t_1t_2\ldots t_i$. Then
$\lambda_i \prec \lambda$.
\item[(iv)] For all $m\in M_x$, we have $l(m) \geq |\phi(m)|$.
\item[(v)] If $\lambda \models x$, we have $l(\delta_\lambda) = |\phi
(\delta_\lambda) |$.
\end{itemize}
\end{lemma}

\begin{proof}
(i) is easy.

(ii) comes from (i) and the following remark: when $z$ and $z'$ are in the same
orbit, multiplying by $t$ splits this orbit into two orbits, thus increasing
by $1$ the number of orbits; when $z$ and $z'$ are in different
orbits, multiplying by $t$ merges their orbits, thus decreasing
by $1$ the number of orbits.

(iii) is an easy induction from (ii) and its proof: as the decomposition
is reduced, the relation $|\sigma| = |t_1t_2\ldots t_{|\sigma|}|$ can only
be achieved if the successive multiplications by the $t_i$ merge orbits.

(iv) and (v) are easy.
\end{proof}

We can now prove proposition \ref{poset}:
\begin{proof}
(i) Let $\lambda,\lambda' \in \P_x$. As the converse implication
is trivial, we only have to check that
$$(\exists m'' \in M_x, \delta_\lambda m''= \delta_{\lambda'})
 \Rightarrow (\exists \lambda'' \in \P_x,
 \delta_{\lambda} \delta_{\lambda''}= \delta_{\lambda'}).$$
Let $m'' \in M_x$ such that $ \delta_\lambda m''= \delta_{\lambda'}$.
Consider the image in $\S_x$ of this identity:
$\phi(\delta_\lambda) \phi(m'')= \phi(\delta_{\lambda'})$.
By lemma \ref{length} (v), we have
$|\phi(\delta_\lambda) | = l(\delta_\lambda)$ and
$|\phi(\delta_{\lambda'}) | = l(\delta_{\lambda'})$,
and, using lemma \ref{length} (iv), 
we have $|\phi(\delta_{\lambda'}) | = 
l(\delta_{\lambda'}) = l( \delta_\lambda m'')
= l( \delta_\lambda) + l(m) \geq |\phi(\delta_\lambda) | + | \phi(m) |$ and
thus $|\phi(\delta_{\lambda'}) | = |\phi(\delta_\lambda) | + | \phi(m) |$. 
Consequently, when concatenating reduced decompositions for 
$\phi(\delta_\lambda)$ and $\phi(m)$, one get a reduced decomposition
for $\phi(\delta_{\lambda'})$. By lemma \ref{length} (iii), this implies 
that the orbit decomposition for $\phi(\delta_{\lambda})$ is finer
than the one for $\phi(\delta_{\lambda'})$, i.e., $\lambda \prec \lambda'$.
We conclude by lemma \ref{diviseurs} (iii).

(ii) By lemma \ref{diviseurs} (iii), the map $D$ is a poset morphism. It
is injective: from $\delta_\lambda$, one recovers $\lambda$ by considering
the orbit decomposition of $\phi(\delta_\lambda)$.
The fact that the inverse map
is a poset morphism, i.e. $\delta_\lambda \prec \delta_{\lambda'} \Rightarrow
\lambda \prec \lambda'$, has already been obtained in our proof of (i).
\end{proof}

\begin{rem} \label{delta} If $x$ is convex, then there is a largest
element $\delta_{\{x\}}$ in the poset $(P_x,\prec)$ (which
corresponds to the coarsest partition which has just one part equal to
$x$).
\end{rem}

\section{Pre-Garside structures and Garside monoids}

The existence of nice normal forms in the Birman-Ko-Lee monoid $B_{\mu_n}$
will be of crucial importance in the proof of our main theorem. This
property is proved in \cite{BKL}, and we could have just quoted and
translated it into our setting. However, it appeared that our 
geometric interpretation allows us to a give a new proof, far less
computational, of some of the main results in \cite{BKL}. 
The Birman-Ko-Lee monoid is a {\em Garside monoid}, in the sense
of \cite{DP}, and this implies (among others) the existence of the normal form. 
To give a simple proof of this fact, we make use a new criterion
of ``Garsiditude''. This criterion relies on the notion of
{\em pre-Garside structure}, which can be seen as an axiomatization of
the context in which most of the proofs in \cite{michel} are actually valid.

When $A$ and $B$ are two sets, we mean by ``partial map'' from $A$ to $B$
a datum consisting of a subset $A' \subset A$ and a map $f:A' \rightarrow B$.
It is convenient to refer to $A'$ implicitly, and to use a slightly abusive
language, e.g. we will write ``$f(a)$ is defined'' instead of ``$a\in A'$''. 

\begin{defi}
Let $P$ be a set.
An atomic partial product on $P$ is 
a partial map $m:P \times P \rightarrow P$ 
(we will denote $m(a,b)$ by $a.b$ or
$ab$),
satisfying the following axioms:
\begin{itemize}
\item[(i)] (Unit element and associativity.)
There exists an element $1 \in P$ such that for all $a\in P$,
 both $1.a$ and $a.1$ exist and are equal to $a$.
For any $a,b,c\in P$, it is equivalent for $ab$ and $(ab)c$ to 
be defined or for $bc$ and $a(bc)$ to be defined
and then $a(bc)=(ab)c$.
\item[(ii)] (Finite number of atoms.)
 Let $P^*= P - \{ 1\}$; the image of $P^*\times P^*$ is in $P^*$,
and the complement $S=P^*-m(P^*\times P^*)$ is finite (the elements
of $S$ are 
the {\em atoms} of $(P,m)$).
\item[(iii)] (Grading.) There exists a function $l:P \rightarrow \NM$ such that
$p\in P^* \Rightarrow l(p)>0$ and $l(ab)=l(a)+l(b)$ whenever $ab$ is defined.
\end{itemize}
\end{defi} 

\begin{defi}
Suppose $P$ is a set together with an atomic partial product.
The associated monoid $M(P)$ is the monoid defined
by the following presentation:
\begin{itemize}
\item As a set of generators we take $P$.
\item For relations we take $ab=c$ whenever $a,b,c\in P$ are such that $ab$ is
defined in $P$ and equal to $c$.
\end{itemize}      
\end{defi}

Note that if $P$ is a subset of a monoid $M'$, and if the partial product
on $P$ is a restriction of the monoid law in $M'$, then there is a natural
morphism $M(P) \rightarrow M'$. If $P$ generates $M'$, the morphism is 
surjective.

As in \cite{michel} we note that we can identify $M(P)$ to the set of 
finite sequences
of elements of $P$, quotiented by the equivalence relation generated by
the equivalence 
when the product $ab$ is defined and equal to $p_i$
of $(p_1,\ldots,p_{i-1},p_i,p_{i+1},\ldots,p_n)$ and
$(p_1,\ldots,p_{i-1},a,b,p_{i+1},\ldots,p_n)$.

Let $(x_1,\ldots,x_n)$ be a sequence of elements of $P$. 
There are $(2(n-1))!/(n! (n-1)!)$ different ways to put brackets
on the product $x_1\ldots x_n$. By an obvious induction from the
associativity axiom, if the product is defined in $P$ for one of these
bracketings, then it is defined in $P$ for any other bracketing, and the 
value of this product does not depend on the choice of the bracketing.
When this is the case, we write $x_1\ldots x_n$ for this product.

\begin{lemma}
Let $(x_1,\ldots,x_n)$ a sequence of elements of $P$ equivalent,
in $M(P)$, to a single term sequence $(y)$. Then 
the product $x_1\ldots x_n$ is defined in $P$ and we have
$x_1\ldots x_n = y$.
\end{lemma}

\begin{proof}
By assumption $(x_1,\ldots,x_n)$ can be transformed into $(y)$ by a
finite rewriting process 
$$l_0= (x_1,\ldots,x_n) \rightarrow l_1 \rightarrow \cdots \rightarrow
l_k= (y)$$ in which, at each step, 
the elementary transformation $l_{j-1} \rightarrow l_{j}$ is 
\begin{itemize}
\item either of the type 
$$(p_1,\ldots,p_{i-1},a,b,p_{i+1},\ldots,p_m) \rightarrow
(p_1,\ldots,p_{i-1},p_i,p_{i+1},\ldots,p_m)$$
\item or of the type
$$(p_1,\ldots,p_{i-1},p_i,p_{i+1},\ldots,p_m) \rightarrow
(p_1,\ldots,p_{i-1},a,b,p_{i+1},\ldots,p_m)$$
\end{itemize}
where $a,b$ are such that their product is defined in $P$ and equal to $p_i$.
Suppose the length $k$ of the rewriting process is minimal.

Suppose  some of the transformations are of the second type, and
chose $j$ maximal such that $l_{j-1} \rightarrow l_{j}$ is of the second
type
$$(p_1,\ldots,p_{i-1},p_i,p_{i+1},\ldots,p_m) \rightarrow
(p_1,\ldots,p_{i-1},a,b,p_{i+1},\ldots,p_m)$$  
As all further steps are of the first type, we have $k=j+m$ and the 
product $p_1\ldots p_{i-1}abp_{i+1}\ldots p_m$ is
defined in $P$ (and is equal to $y$).  
Choosing a bracketing starting by  $\ldots(ab)\ldots$, we see that
the product $p_1\ldots p_{i-1}p_ip_{i+1}\ldots p_m$ must also be defined
in $P$ (and equal to $y$). But this yields a rewriting of length
$j+m-2 < k$ and we have a contradiction.

Thus there are no transformation of the second type in a minimal rewriting.
The result follows.
\end{proof}

The following propostion is a straightforward consequence of the lemma:

\begin{prop}
\label{divides}
\begin{itemize}
\item[(i)] The natural map from $P$ to $M(P)$ is injective.
\item[(ii)] If $x\in M(P)$ divides $a\in P$ then $x\in P$.
\end{itemize}
\end{prop}    

%

\begin{defi}
A {\em pre-Garside} structure on a set $P$ is an atomic partial product,
with set of atoms $S$, satisfying the following additional axioms:
\begin{itemize}
\item[(iv)] If two elements of $S$ have a common right
multiple in $P$, they have
a least common right multiple. When $s,t\in S$ have a least right common 
multiple, we write it $\Delta_{s,t}$.
\item[(iv')] If two elements of $S$ have a common left
multiple in $P$, they have
a least common left multiple. 
\item[(v)] If $s,t \in S$ have a common right multiple in $P$, and if $a\in P$
is such that $as\in P$ and $at \in P$, then $a\Delta_{s,t}\in P$.
\item[(vi)] For all $m\in M(P)$ and $a,b\in P$, if either $am=bm$ or
$ma=mb$, then $a=b$.
\end{itemize}
\end{defi}


In \cite{michel}, $S$ was taken to be the set of
usual Artin generators, and $P$
was the set of reduced braids. To handle the Birman-Ko-Lee monoid, we will
take for $S$ the set $S_{\mu_n}$ and for $P$ the set $P_{\mu_n}$.

Condition (v) is satisfied {\it e.g.} if there exists a common right
multiple in $P$ of all elements of $P$, which is the case in the usual braid
monoid and also in the Birman-Ko-Lee
monoid (see \ref{delta}).
Note that such an element is necessarily unique and if it exists, 
the finiteness of $S$ implies that $P$ is finite. We shall see
(\ref{Delta existe}) that
its existence is equivalent to the fact that all pairs of elements of $S$ have
a right lcm.


The following lemmas and propositions are rephrasings of
\cite{michel}, 1.4 to 1.9.
\begin{lemma}
\label{1.4}
Let $X$ be a finite subset of $M(P)$ such that
\begin{itemize}
\item If $x\in X$, $a\in M(P)$, $a\prec x$ then $a\in X$.
\item If $a\in M(P)$, $s,t\in S$, $as, at\in X$, then $\Delta_{s,t}$
exists and $a\Delta_{s,t}\in X$.
\end{itemize}
Then there exists $g\in X$ such that $X$ is the set of left divisors of $g$.
\end{lemma}
\begin{proof}
The statement follows \cite{michel} 1.4, and the proof is exactly the same.
\end{proof}
\begin{prop}
\label{pgcd dans P}
Any two elements of $P$ have a left g.c.d.\ in $P$.
\end{prop}
\begin{proof}
We can follow the proof of \cite{michel} 1.6, replacing 1.5 in {\it
loc.\ cit.} by axiom (iv).
\end{proof}
We shall denote by $a\wedge b$ the left g.c.d.\ of $a$ and $b$.
\begin{prop}
\label{alpha dans P}
For $a$ and $b$ in $P$ there exists a unique maximal $c\prec b$ such that
$ac\in P$. 
\end{prop}
\begin{proof} 
We apply
lemma \ref{1.4} to the set $X$ of $c$ such that $c\prec b$ and $ac\in P$.
To check the
assumptions of that lemma, we need that if $s,t\in S$ and if $cs,ct\in
X$ then $\Delta_{s,t}$ exists and
$ac\Delta_{s,t}\in P$ (then we have $c\Delta_{s,t}\in X$).
Since $cs$ and $ct$ divide $b\in P$,
by cancellability and \ref{divides} $s$ and $t$ have a common multiple in $P$,
so by axiom (iv) $\Delta_{s,t}$ exists and by axiom (v)
$ac\Delta_{s,t}\in P$.
\end{proof}
\begin{defi}
\label{alpha2 et omega2}
In the situation of \ref{alpha dans P} we denote $\alpha_2(a,b)$ the element
$ac$ and we denote $\omega_2(a,b)$ the unique $d\in P$
such that $b=cd$. We thus have $ab=\alpha_2(a,b)\omega_2(a,b)$.
\end{defi}
Note that the uniqueness of $\omega_2$ follows from axiom (vi).
\begin{prop}
\label{alpha2(ab,c)}
For $a,b,c,ab\in P$ we have $\alpha_2(ab,c)=\alpha_2(a,\alpha_2(b,c))$.
\end{prop}
\begin{proof}
The statement is \cite{michel}, 1.8 and the proof is the same.
\end{proof}
\begin{prop}
\label{omega2(ab,c)}
For $a,b,c,ab\in P$ we have 
$$\omega_2(ab,c)=\omega_2(a,\alpha_2(b,c))\omega_2(b,c).$$
\end{prop}
\begin{proof}
By propositions \ref{alpha dans P}
and \ref{alpha2(ab,c)} the products of both sides with
$\alpha_2(ab,c)$ are equal. By axiom (vi) we will be done if we show
that $\omega_2(a,\alpha_2(b,c))\omega_2(b,c)$ is in $P$.
By definition \ref{alpha2 et omega2} there exists $z\in P$ such that
$\alpha_2(b,c)=bz$ and $c=z\omega_2(b,c)$. As $ab\prec abz$ and $ab\in P$ we have
$\alpha_2(a,bz)=abz_1$ for some $z_1\prec z$ (by definition of $\alpha_2$ and
cancellability in $P$). Hence $bz=bz_1\omega_2(a,bz)$ and
$z=z_1\omega_2(a,bz)$. So
$z_1\omega_2(a,\alpha_2(b,c))\omega_2(b,c)=z_1\omega_2(a,bz)\omega_2(b,c)=
z\omega_2(b,c)=c$. The result follows since any divisor of $c$ is in $P$ by
\ref{divides}.
\end{proof}

We now extend the definition of $\alpha_2$ to $M(P)$, following
\cite{michel}, 2.1 to 2.6. All the proofs of {\it loc.\ cit.} can be reproduced,
replacing proposition 1.5 in {\it loc.\ cit.} by axiom (iv).
\begin{prop}
\label{alpha}
There is a unique function $\alpha: M(P)\to P$ extending the identity of $P$
and satisfying $\alpha(ab)=\alpha_2(a,b)$ for $a,b\in P$ and
$\alpha(gh)=\alpha(g\alpha(h))$ for $g,h\in M(P)$. Moreover $\alpha(g)$ is the
unique maximal element in $\{c\in P| c\prec g\}$.
\end{prop}
\begin{prop}
\label{omega}
There exists a unique function $\omega: M(P)\to M(P)$ such that
$\omega(ab)=\omega(a,b)$ for $a,b\in P$ (in particular
$\omega(a)=1$ for $a\in P$) and such that
$\omega(gh)=\omega(g\alpha(h))\omega(h)$ for $g,h\in M(P)$.
\end{prop}
\begin{prop}
\label{2.3}
Let $g\in M(P)$; then $\omega(g)$ is the unique $y\in M(P)$ such that
$g=\alpha(g)y$.
\end{prop}
\begin{prop}
\label{2.4}
The monoid $M(P)$ has left and right cancellation property (i.e., in axiom
(vi) we can replace the condition $a,b\in P$ by $a,b\in M(P)$).
\end{prop}
\begin{prop}
\label{2.5}
If $s,t\in S$ divide $a\in M(P)$ on the left, then $\Delta_{s,t}$ exists and
divides $a$.
\end{prop}
\begin{prop}
\label{2.6}
For $a,b\in M(P)$ there exists a unique maximal $c\in M(P)$ (for $\prec$)
such that $c\prec a$ and $c\prec b$.
\end{prop}
We shall still denote this left g.c.d.\ by $a\wedge b$.
\begin{prop}
\label{lcm in M}
A family of elements
of $M(P)$ which has a right (resp.\ left) common multiple 
has a right (resp.\ left) lcm in $M(P)$.
\end{prop}
\begin{proof}
Assume that all elements of the family $\{a_i\}_{i\in I}$ divide $c$.
If we can
apply \ref{1.4} to the set $X$ of
elements of $M(P)$ which divide all common multiples of the $a_i$, it will
give the result. Let us check the assumption of \ref{1.4}.
This set $X$
is finite as it is included in the set of divisors of $c$. The first
assumption of \ref{1.4} is clearly satisfied. The second assumption is a
consequence of the fact that
if $s,t\in S$ are such that
$xs$ and $xt$ divide some element $xz$, then $s$ and
$t$ divide $\alpha(z)$, so $\Delta_{s,t}$ exists and divides $z$, whence
$x\Delta_{s,t}$ divides $xz$.
\end{proof}
We have a more precise result for elements of $P$:
\begin{prop}
\label{lcm in P}
If a family of elements of $P$ has a common right (resp.\ left)
multiple in $M(P)$ then
its right (resp.\ left) lcm exists and is in $P$.
\end{prop}
\begin{proof}
The lcm exists by the preceding proposition.
Let $m$ be this lcm. Any divisor of $m$ divides
$\alpha(m)\in P$, whence the result.
\end{proof}

We now get a normal form for any element of $M(P)$ exactly as in
\cite{michel}.
\begin{defi} 
\label{normal form}
A decomposition $(g_1,\ldots,g_n)$ of an element $g_1\ldots g_n$ of $M(P)$ is
said to be its normal form if no $g_i$ is equal to 1 and for any $i$
we have $g_i=\alpha(g_i\ldots g_n)$.
\end{defi}
The following statement is 4.2 of \cite{michel} and the same proof applies.
\begin{prop}
\label{local}
A decomposition $(g_1,\ldots ,g_i)$ with $g_i\in P$ is a normal form
if and only if $(g_i,g_{i+1})$ is a normal form for any $i$.
In particular any segment
$(g_i,\ldots,g_j)$ of a normal form $(g_1,\ldots,g_n)$ is a normal form.
\end{prop}
In the same way, statements 4.5 to 4.9 and 5.1 to 5.3 of \cite{michel}
generalize to our setting.

Assume that
all elements of $P$ have a
right common multiple (which is in $P$ by \ref{lcm in P}).
By \ref{lcm in M}, this is the same as assuming that any
finite subset of elements of $M(P)$ has a lcm.
As already noticed, the existence of a lcm $\Delta$ of all elements in
$P$ implies the
finiteness of $P$. The converse is not true, but we have:
\begin{prop}
\label{Delta existe}
The elements of $P$ have a common right multiple if and only if $P$ is finite
and any pair of elements of $S$ has a common multiple.
\end{prop}
\begin{proof}Implication ``only if'' is clear. The converse is
an immediate application of \ref{1.4}, taking $X=P$.
\end{proof}
By \ref{alpha}, for any $a\in M(P)$ we have $\alpha(a)=a\wedge\Delta$.
\begin{prop}
\label{automorphisme} Assume $P$ has a 
 right lcm $\Delta$.
\begin{itemize}
\item
There is an automorphism $x\mapsto\bar x$ of $M(P)$ mapping $S$ to itself
and such that $x\Delta=\Delta\bar x$. 
\item 
The element $\Delta$ is the left lcm\ of $P$.
\end{itemize}
\end{prop}
\begin{proof} It is sufficient to define the automorphism on $P$.
As $\Delta$ is a right multiple of all elements of $P$,
for any $a\in P$ there exists a unique
$\Delta_a\in P$ such that $\Delta=a\Delta_a$ and
there exists a unique 
$\bar a\in P$ such that $\Delta=\Delta_a\bar a$, so that
$a\Delta=a\Delta_a\bar a=\Delta\bar a$.
The map $a\mapsto \bar a$
is injective by cancellability in $M(P)$ and is thus compatible with the
product. As $P$ is finite, it is surjective.
If $a\in S$, by surjectivity $\bar a$ cannot be the product of two non
trivial elements of $P$, so has to be in $S$ by axiom (ii).

As $a\mapsto \bar a$ is a bijection of $P$,
the above proof shows that $\Delta$
is a left multiple of all elements of $P$.
So it is the left lcm\ of $P$.
\end{proof}

We can now compare our formalism with the one in \cite{DP}:
if $P$ is a pre-Garside structure with right common multiples, then, it is
readily seen that any two elements of $M(P)$ have left and 
right common multiples:
indeed, if $m_1,m_2\in M(P)$ are both product of $n$ or less elements of
$S$, then $\Delta^n$ is a left and 
right multiple of both $m_1$ and $m_2$ (use 
proposition \ref{Delta existe}:
the conjugation by $\Delta$ maps $S$ into itself).
Thus, by proposition \ref{lcm in M}, pairs of elements in $M(P)$ have both
a right lcm and a left lcm. Moreover, by proposition \ref{Delta existe},
$\Delta$ is both the right lcm and the left lcm of $P$. Conversely, elements
of $P$ are left and right divisors of $\Delta$, and, by proposition
 \ref{divides}, $\Delta$ has no other left or right divisors. 
We have proved that $M(P)$ is a Garside monoid, as defined in \cite{DP}. 

Conversely, if $M$ is a Garside monoid, then the restriction of the monoid
product to the set $P$ of {\em simple 
elements} (see \cite{DP}) is a pre-Garside structure with common multiples.

We have proved the:

\begin{theo} Let $P$ be set, endowed
with a pre-Garside structure, such that all atoms have a common right multiple.
Then $M(P)$ is a Garside monoid.

Conversely, when $M$ is a Garside monoid with
fundamental element $\Delta$, the set $P$ of divisors of $\Delta$
has a pre-Garside structure for the partial product obtained by restriction
of the product in $M$; all elements of $P$ have a common right multiple and
$M\simeq M(P)$.
\end{theo}

Note that our approach already gives interesting results (e.g. \ref{2.4}
and \ref{points fixes} below) even when some elements of $P$ have
no common multiple, e.g. in the case of braid groups attached to infinite
type Coxeter groups; in this case, $M(P)$ is no longer a Garside monoid.
More important is, for the present article, that pre-Garside structures
provide us with a convenient criterion to check that the Birman-Ko-Lee monoid
is a Garside monoid.

Let us recall the following proposition
about Garside monoids, which is already in \cite{DP}.


\begin{prop}
\label{injection du monoide}
Assume that $P$ has a right lcm.
The monoid $M(P)$ injects into the group $G(P)$
having same presentation and any element of $G$ can be written uniquely 
$x^{-1} y$ with $x,y\in M(P)$ and $x\wedge y=1$.
\end{prop}
\begin{proof}
The statement follows \cite{michel} 3.2 and the same proof applies.
\end{proof}
We now generalize \cite{michel} 4.4.
\begin{prop}
\label{points fixes} 
Let $\Gamma$ be
a group of automorphisms of $M(P)$ stabilizing $S$. Let $\Sigma$
be the set of lcms of $\Gamma$-orbits in $S$ which exist and
are not the product of other such lcms;
then $P^\Gamma$ has a pre-Garside structure with atoms $\Sigma$; the monoid
$M(P)^\Gamma$ identifies with $M(P^\Gamma)$.
If moreover we assume that $P$ has a right lcm
then the group $G(P)^\Gamma$
identifies with $G(P^\Gamma)$ ({\it cf.} \ref{injection du monoide}).
\end{prop}
\begin{proof} Axiom (i) holds as it holds in $P$ and as the product of two
$\Gamma$-fixed elements is $\Gamma$-fixed.
Let $x\in P^\Gamma$, and let $s\in S$ such that $s\prec x$;
then $x$ is divisible
by all elements in the orbit of $s$, so is divisible by their lcm
(which exists).
So $P^\Gamma$ is generated by the lcms of $\Gamma$-orbits in $S$ which
exist, thus
by $\Sigma$, and we have arranged for elements of $\Sigma$ to be atoms,
so axiom (ii) is satisfied.
The length inherited from $P$ is
still compatible with the product so we have axiom (iii) (but note that the
elements of $\Sigma$ may have length greater than 1 even if all elements of
$S$ had length 1).
The lcm of two elements of $P^\Gamma$ is
$\Gamma$-fixed by its uniqueness and is in $P$ by \ref{lcm in P}, whence 
axioms (iv) and  (iv'). 
If $\sigma$ and $\tau$ are in $\Sigma$ and have a lcm
$\Delta_{\sigma,\tau}$
and if $a\in P^\Gamma$ is such that $a\sigma$ and $a\tau$ are in
$P$ then $a\Delta_{\sigma,\tau}$ is the lcm of $a\sigma$ and $a\tau$
so is in $P$ by \ref{lcm in P}, and is clearly $\Gamma$-fixed,
whence axiom (v).
The monoid $M(P^\Gamma)$ is by definition a submonoid of $M(P)$,
so axiom (vi) holds.

Let $x\in M(P)^\Gamma$;
the uniqueness of its normal form implies that each term
is in $P^\Gamma$. On the other hand, if $x\in M(P)^\Gamma$,
as $\alpha(x)$ is the unique maximal element in $P$ dividing $x$ it is also
the unique maximal element in $P^\Gamma$ dividing $x$, so, by
the definition of normal forms, the normal form of an element of
$M(P)^\Gamma$ is a normal form in $M(P^\Gamma)$.
This shows that $M(P)^\Gamma$ identifies with $M(P^\Gamma)$. The same argument
shows that $G(P)^\Gamma$ identifies with $G(P^\Gamma)$ when $\Delta$
exists, as $\Delta$ is in $P^\Gamma$.
\end{proof}

\section{The Birman-Ko-Lee monoid}

Let $x\in X_n$, and let $S_x$, $P_x$ and $M_x$ as defined as in section
\ref{localmonoid}.

\begin{lemma}
The restriction of the product in $M_{x}$ defines an atomic partial product
in $P_{x}$ with set of atoms $S_x$.
\end{lemma}

\begin{proof}
Axiom (i) is a consequence of the associativity in the monoid $M_{x}$.
Axioms (ii) and (iii) are easy consequences of the existence of the
length function $l: M_x \rightarrow \NM$ introduced in section
\ref{localmonoid}. The atoms are clearly, by construction,
 the elements  of $S_x$.      
\end{proof}

We call
{\em Birman-Ko-Lee monoid} the monoid $M_{\mu_n}$ (whenever $x$ is convex,
we have $M_x \simeq M_{\mu_n}$).

{\bf \flushleft Remark.}
We can map the base point $\mu_n$ to the usual base point $\{1,\ldots,n\}$ by
choosing $z_1\in\mu_n$ and numbering clockwise the elements
of $\mu_n$ starting at $z_1$ as $z_1,z_2,\ldots,z_n$;
then the $\de(z_i,z_{i+1})$ map to the usual Artin generators
$\sigma_i$ of $\pi_1(X_n,\mu_n)$  and the 
$\de(z_i,z_j)$ map to the generators $a_{ij}$ considered in \cite{BKL} (note
that they always suppose $i>j$ while we don't order the pairs $\{z,z'\}$).
Then the elements of $P_{\mu_n}$ map to the {\em canonical factors}
defined in \cite{BKL}. In order to compare
our definition to that of \cite{BKL}
one should note that \cite{BKL} call $a_{ij},a_{kl}$ an ``obstructing pair''
exactly when the partition $\{\{z_i,z_j\},\{z_k,z_l\}\}$ is obstructing.

The monoid $M_{\mu_n}$ enjoys remarkable properties which are consequence of
the discussion in the previous section and of the following theorem.
\begin{theo}
The restriction of the product in $M_{\mu_n}$ to a partial product
in $P_{\mu_n}$ is a pre-Garside structure, and $P_{\mu_n}$ has a
lcm $\delta$.
\end{theo}

\begin{proof}
Since $\mu_n$ is convex, $\delta=\delta_{ \{\mu_n \} }$ is an lcm of
$P_{\mu_n}$, as explained in \ref{delta}.
This eliminates the need to check axiom (v).

Let us now prove axiom (iv). Using the isomorphism in proposition \ref{poset},
and noticing that any pair of elements of $\mu_n$ defines a non obstructing
edge, 
we have to prove:
for all $z_1,z_1',z_2,z_2' \in \mu_n$, the set 
$$E = \{ \lambda \models \mu_n | \{\{ z_1,z_1'\} \} \prec \lambda 
\; \text{and} \; \{\{ z_2,z_2'\} \} \prec \lambda \}$$
has a minimum element for $\prec$.
We discuss by cases:
\begin{itemize}
\item Either $[z_1,z_1'] \cap [z_2,z_2'] = \emptyset$. Then let
$\lambda = \{\{ z_1,z_1'\}, \{ z_2,z_2'\} \}$. Clearly, $\lambda \in E$,
and any partition coarser than $\{\{ z_1,z_1'\} \}$ and $ \{\{ z_2,z_2'\} \}$
must be coarser than $\lambda$.
\item  Either $[z_1,z_1'] \cap [z_2,z_2'] \neq \emptyset$. Let
$\lambda = \{\{ z_1,z_1', z_2,z_2'\} \}$ (the non-trivial part may have
three or four elements). Clearly $\lambda \in E$. Now let $\lambda' \in E$.
Consider $\nu_1$ the part of $\lambda'$ in which $z_1$ lies,
$\nu_2$ the part of $\lambda'$ in which $z_2$ lies. Because 
$\{\{ z_1,z_1'\} \} \prec \lambda'$
and $\{\{ z_2,z_2'\} \} \prec \lambda'$, we have $z_1'\in \nu_1$ and
$z_2'\in\nu_2$. As $\lambda'$ is non obstructing, we must have 
$\nu_1=\nu_2$ (otherwise we would have $\overline{\nu_1}\cap
\overline{\nu_2} = \emptyset$, which contradicts 
$[z_1,z_1'] \cap [z_2,z_2'] \neq \emptyset$). 
Thus $\lambda$ is finer than $\lambda'$.
\end{itemize}

(iv') is proved similarly.
 
To prove (vi), consider the natural morphism
 $M(P_{\mu_n}) \rightarrow M_{\mu_n}$ and its composition with
the epimorphism $M_{\mu_n} \rightarrow \S_{\mu_n}$.
Let $m\in M(P_{\mu_n})$, $\lambda,\lambda'\models \mu_n$ such that
either $\delta_\lambda m = \delta_{\lambda'} m $ or
$m \delta_\lambda  = m \delta_{\lambda'} $. Denote by $\sigma$ the image
of $m$ in $\S_{\mu_n}$. We have
$\phi(\delta_\lambda) \sigma = \phi(\delta_{\lambda'}) \sigma$ or
$\sigma \phi(\delta_\lambda) = \sigma \phi(\delta_{\lambda'})$. As
$\S_{\mu_n}$ is a group, this implies
$\phi(\delta_\lambda) = \phi(\delta_{\lambda'})$. By lemma \ref{phiinjective},
the restriction of $\phi$ to $P_{\mu_n}$ is injective, thus we have as
required $\delta_\lambda = \delta_{\lambda'}$.
\end{proof}

{\bf \flushleft Remark.} According to the isomorphism in proposition
 \ref{poset},
the above theorem implies that the poset of non obstructing partitions
of $\mu_n$ is a lattice. If $\lambda_1,\lambda_2\models \mu_n$, 
the set $$E=\{\lambda\models\mu_n | \lambda_1\prec \lambda \; \text{and}
\; \lambda_2\prec\lambda\}$$ admits a minimum element. We leave to the reader
to check that this can be proved directly, using the following arguments:
given two (possibly obstructing) partitions, the set of 
(possibly obstructing) partitions
coarser than the two partitions admits a minimum
element; given a (possibly obstructing) partition, the set of 
non obstructing partitions which are coarser admits a minimum element. 

Of course, by construction, $M_{x}$ embeds in $B_{x}$ for all $x\in X_n$.
What is specific to the Birman-Ko-Lee monoid is the following result,
which is the analog in our intrinsic setting of the embedding theorem in
\cite{BKL}. 

\begin{coro}
The natural morphism $M(P_{\mu_n}) \rightarrow M_{\mu_n}$ is an isomorphism.
\end{coro}

\begin{proof}
Consider the group $G(P_{\mu_n})$ given by the group presentation corresponding
to the monoid presentation of $M(P_{\mu_n})$. As the partial product in 
$P_{\mu_n}$ is a restriction of the one in $B_{\mu_n}$,
 there is a natural morphism
$$G(P_{\mu_n}) \rightarrow B_{\mu_n}.$$
By proposition \ref{Sgenerates}, this morphism is surjective. 

Note that it is easy to adapt Sergiescu's 
presentations to find a presentation of $B_{\mu_n}$ where the generators
corresponds to the elements of $S_{\mu_n}$. A set of defining relations
is for example given in \cite{BKL}, proposition 2.1.:
\begin{itemize}
\item $\delta_{\{ z_1,z_1' \}  } \delta_{\{ z_2,z_2'\} } =
   \delta_{\{ z_2,z_2' \}  } \delta_{\{ z_1,z_1'\} }$ when
$\{ \{ z_1,z_1' \} , \{ z_2,z_2'\} \}$ is non obstructing,
\item $\delta_{\{ z_1,z_2 \}  } \delta_{\{ z_2,z_3\} }=
\delta_{\{ z_2,z_3 \}  } \delta_{\{ z_3,z_1\} }=
\delta_{\{ z_3,z_1 \}  } \delta_{\{ z_1,z_2\} }$ for $z_1,z_2,z_3$ coming
in clockwise order.
\end{itemize}
These relations are valid in $P_{\mu_n}$, thus in $G(P_{\mu_n})$, and the
morphism $ G(P_{\mu_n}) \rightarrow B_{\mu_n} $ is an isomorphism.

We conclude using the natural commutative diagram:
$$\xymatrix{ M(P_{\mu_n}) \ar@{>>}[r] \ar[d] & M_{\mu_n} \ar@{^{(}->}[d]\\
          G(P_{\mu_n}) \ar[r]^\sim & B_{\mu_n} }$$
and the injectivity of the map $M(P_{\mu_n}) \rightarrow  G(P_{\mu_n}) $
(proposition \ref{injection du monoide}).
\end{proof}

We identify $M(P_{\mu_n})$ and $M_{\mu_n}$ through the natural isomorphism.
As $M_{\mu_n}$ is a Garside monoid, we have nice normal forms, an algorithm
for the word problem, $\ldots$ The following property is the one we need for
the proof of our main theorem.

\begin{prop}\label{conj delta}
Let $\zeta=e^{\frac{2i\pi}n}$.
\begin{itemize}
\item[(i)] For $\lambda\models\mu_n$ the automorphism of $M_{\mu_n}$ induced by
$\delta$ maps $\delta_\lambda $ to $\delta_{\zeta\lambda}$.
\item[(ii)] For $d$ dividing $n$,
the centralizer of $\delta^{n/d}$ in $\pi_1(X_n,\mu_n)$ is
generated by the elements $\delta_\lambda$ for $\lambda\models\mu_n$ such that
$\zeta^{n/d}\lambda=\lambda$.
\end{itemize}
\end{prop}
\begin{proof}
It is an immediate consequence of \ref{tourner} that (i) holds for a generator
of the form  $\de(s_i,s_{i+1})$ where $s_i$ and $s_{i+1}$ are two consecutive
points in a numbering of $\mu_n$. It follows that it holds for any generator
$\de(s_1,s_a)$ by using that
$\de(s_1,s_a)=\de(s_1,s_2)\ldots\de(s_{a-1},s_a)\de(s_{a-2},s_{a-1})^{-1}\ldots
\de(s_1,s_2)^{-1}$ (which follows also from \ref{tourner}), and it follows thus
for any element of $P$.
Part (ii) is then a consequence of \ref{points fixes}.
\end{proof}

{\bf \flushleft Remark.}
Note that \ref{points fixes} gives a  ``Birman-Ko-Lee'' presentation of
$C_\bB(\delta^i)$. Let us work out an example to show that one has to
take only a part of the lcms of $\delta^i$-orbits on $S$ to get an
atomic set. Take $n=6$ and $i=2$. Then the $\delta^2$-orbits in $S$
are $$\displaylines{\{\de(s_1,s_2),\de(s_3,s_4),\de(s_5,s_6)\},
\{\de(s_2,s_3),\de(s_4,s_5),\de(s_1,s_6)\},\hfill\cr
\hfill\{\de(s_1,s_3),\de(s_3,s_5),\de(s_1,s_5)\}\hbox{ and }
\{\de(s_1,s_4),\de(s_2,s_5),\de(s_3,s_6)\}\cr}$$ whose respective lcm are
$$\delta_{(\{s_1,s_2\},\{s_3,s_4\},\{s_5,s_6\})},
\delta_{(\{s_2,s_3\},\{s_4,s_5\},\{s_1,s_6\})},
\delta_{\{s_1,s_3,s_5\}}\hbox{ and }\delta;$$ but
$\delta=\delta_{\{s_1,s_3,s_5\}}\delta_{(\{s_1,s_2\},\{s_3,s_4\},\{s_5,s_6\})}$,
so must be eliminated.

\section{A geometric normal form for canonical factors.}

Let $n$ be a positive integer. If $K_1,K_2$ are non-intersecting 
compact subsets
of $\CM$, we denote by $d(K_1,K_2)$ the positive number
$$\inf_{(z_1,z_2)\in K_1\times K_2} |z_1-z_2|.$$
It is clear that 
$$\inf_{\lambda\models\mu_n} \inf_{\stackrel{\nu_1,\nu_2\in \lambda}
{\nu_1\neq\nu_2}
} d(\overline{\nu_1},\overline{\nu_2}) > 0 .$$

We define
$$\epsilon_n := \frac{1}{3} \inf_{\lambda\models\mu_n} \inf_{\stackrel{\nu_1,\nu_2\in \lambda}
{\nu_1\neq\nu_2}
} d(\overline{\nu_1},\overline{\nu_2}).$$
The exact value is not important, we will only use the fact that
$\epsilon_n$ is fixed and small enough. In particular, the reader should check
for himself that in the following definition,
the strings do not intersect and thus $\gamma_\lambda$ is well defined
(note that the function $f$ has been defined in the first section).

\begin{defi}
Let $\lambda\models\mu_n$. 
\begin{itemize}
\item Let $z\in \mu_n$.
We define a path $\gamma_z:[0,1] \rightarrow \CM$ in the
following way:
\begin{itemize}
\item If $z$ is not in the support of $\lambda$, then we set
$$\forall t, \gamma_z(t)=z.$$
\item If $z$ is in a part of $\lambda$ with exactly two elements $z$ and
$z'$ we set:
$$\forall t, \gamma_z(t)= f(z,z',\epsilon_n)(t)$$
\item If $z$ is in a part $\nu$ of $\lambda$ with three or more elements,
we denote by $z'$ the element
of $\nu$ immediately after $z$ in the direct (i.e., anti-clockwise)
 cyclic order on the
vertices of $\overline{\nu}$ and we set:
$$\forall t, \gamma_z(t)=z+ t(z'-z).$$
\end{itemize}
\item These strings uniquely determine an element of $\Omega(X_n,\mu_n)$
which we denote by $\gamma_\lambda$. 
\end{itemize}
\end{defi}

\begin{lemma}
\label{gamma}
Let $\lambda\models\mu_n$, and let $\zeta=e^{\frac{2i\pi}{n}}$.
\begin{itemize}
\item[(i)] The loop $\gamma_\lambda$ represents
$\delta_\lambda$.
\item[(ii)]
 $\forall z\in \mu_n, \gamma_{z\cdot\lambda}=z\cdot\gamma_\lambda$.
\item[(iii)] Let $\lambda'\models\mu_n$. We have
$\lambda = \lambda' \Leftrightarrow \delta_\lambda = \delta_{\lambda'} 
\Leftrightarrow \gamma_\lambda = \gamma_{\lambda '}$.
\item[(iv)] Let $k\in \NM$. Denote by $d$ the order of the root of unity
$\zeta^k$. Then
$$\zeta^k \lambda = \lambda 
\Leftrightarrow 
\delta^k \delta_\lambda = \delta_\lambda \delta^k \Leftrightarrow 
\gamma_\lambda \in \Omega(X_n^{\mu_d},\mu_n).$$
\end{itemize}
\end{lemma}

\begin{proof}
The first three statements are obvious consequences of the previous 
definition.

As $\delta^{-k} \delta_\lambda \delta^{k}= \delta_{\zeta^k\lambda}$,
we have, using (iii),
$$e^{\frac{2ik\pi}{n}} \lambda = \lambda
\Leftrightarrow
\delta^k \delta_\lambda = \delta_\lambda \delta^k \Leftrightarrow
\gamma_{\zeta^k \lambda} = \gamma_{\lambda}.$$
By (ii), we have $$\gamma_{\zeta^k \lambda} = \gamma_{\lambda}
\Leftrightarrow \zeta^k \gamma_{\lambda} = \gamma_{\lambda}.$$ It 
is clear that $\zeta^k \gamma_{\lambda} = \gamma_{\lambda}$ is equivalent
to $\gamma_\lambda \in \Omega(X_n^{\mu_d},\mu_n)$. We have proved (iv).
\end{proof}

\section{Proof of the theorem}

We start with one particular case.

\begin{prop}
\label{case1}
Let $r,d\in\NM$. Let $n=dr$.
The inclusion $X_n^{\mu_d} \subset X_n$ induces an isomorphism
$$\pi_1(X_n^{\mu_d},\mu_n) \stackrel{\sim}{\longrightarrow} 
\pi_1(X_n,\mu_n)^{\mu_d}.$$
\end{prop}

\begin{proof}
The injectivity part of the proposition has been proved in \cite{bessis}.
Let us prove the surjectivity. 
Let $\zeta=e^{\frac{2i\pi}{d}}$.
Multiplication by $\zeta$ coincides with
conjugating by $\delta^{n/d}$, so,
by \ref{conj delta}(ii), we know that $\pi_1(X_n,\mu_n)^{\mu_d}$ is generated by
the $\delta_{\lambda}$ such that $\zeta\lambda=\lambda$.
Using (iv) and (v) of lemma \ref{gamma}, this means that
$\gamma_\lambda\in \Omega(X_n^{\mu_d},\mu_n)$.
Thus any such $\delta_\lambda$ is in the image of
$\pi_1(X_n^{\mu_d},\mu_n)$. So
we have proved that the image of $\pi_1(X_n^{\mu_d},\mu_n)$ is
$\pi_1(X_n,\mu_n)^{\mu_d}$.
\end{proof}

We now have to study the other type of regular numbers for $\S_n$.
We denote by $\nu_n$ the basepoint of $X_n$ defined by
$\nu_n:=\mu_{n-1}\cup\{ 0 \}$.

\begin{prop}
Let $r,d\in\NM$. Let $n=dr+1$.
The inclusion $X_n^{\mu_d} \subset X_n$ induces an isomorphism
$$\pi_1(X_n^{\mu_d},\nu_n ) \stackrel{\sim}{\longrightarrow}
\pi_1(X_n,\nu_n)^{\mu_d}.$$
\end{prop}

\begin{proof}

For all integer $m$, let us denote by $X^*_m$ the space of configurations
of $m$ points in $\CM^*$, with its natural topology.

There is a natural inclusion $X^*_m \subset X_m$ and a natural injection
$X^*_m \subset X_{m+1}$ defined by $x\mapsto x\cup \{ 0\}$.
The action of $\mu_d$ on $X_m$ restricts to an action on $X^*_m$.
Thus we have a commutative diagram of continuous maps:

$$\xymatrix{ X_n^{\mu_d} \ar@{^{(}->}[d] & X_{dr}^{*\mu_d} \ar[l] \ar[r]
\ar@{^{(}->}[d] & X_{dr}^{\mu_d} \ar@{^{(}->}[d] \\
X_n & X_{dr}^* \ar[r] \ar[l] & X_{dr} } $$

According to \cite{bessis}, lemme 3.1., the first line consists of
homeomorphisms. Consider the following commutative diagram of group
morphisms:
$$\xymatrix{ \pi_1(X_n^{\mu_d},\nu_n)
 \ar[d] \ar@/_3em/[dd]_\alpha
 & \pi_1(X_{dr}^{*\mu_d},\mu_{dr}) \ar[l]_\sim \ar[r]^\sim
\ar[d] & \pi_1(X_{dr}^{\mu_d},\mu_{dr}) \ar[d] \ar@/^3em/[dd]^\beta \\
 \pi_1(X_n,\nu_n) & \pi_1(X_{dr}^*,\mu_{dr}) \ar[r]_B \ar[l]^A & 
\pi_1(X_{dr},\mu_{dr})  \\
 \pi_1(X_n,\nu_n)^{\mu_d} \ar@{^{(}->}[u] &
 \pi_1(X_{dr}^*,\mu_{dr})^{\mu_d} \ar@{^{(}->}[u] 
 \ar[r]_{b} \ar[l]^{a} &
\pi_1(X_{dr},\mu_{dr})^{\mu_d} \ar@{^{(}->}[u] 
 }, $$
where $\alpha$ and $\beta$ are defined by functoriality, as
in the introduction. According
to our proposition \ref{case1}, $\beta$ is an isomorphism.
What we want to prove is that $\alpha$ is also
an isomorphism. This will result from the fact that both $a$ and
$b$ are isomorphisms.
\begin{itemize}
\item The map $a$ is an isomorphism: by an easy (and almost
classical) argument, one can see that $A$ is injective and identifies
$\pi_1(X_{dr}^*,\mu_{dr})$ with the subgroup of $\pi_1(X_n,\nu_n)$ consisting
of elements whose associated permutations of $\nu_n$ fix the 
point $0$. As clearly $\pi_1(X_n,\nu_n)^{\mu_d}$ is included in this
subgroup, $a$ is an isomorphism.
\item The map $b$ is an isomorphism: as noticed in \cite{bessis} (fact
$(*)$, used in the proof of theorem 3.2.(I), page 14), 
there is an exact sequence 
$$\xymatrix@1{ 1 \ar[r] & F_n \ar[r] & \pi_1(X_{dr}^*,\mu_{dr}) \ar[r]_B &
\pi_1(X_{dr},\mu_{dr}) \ar[r] & 1 }$$
where $F_n$ is the free group on $n$ generators and such that the action
of $\mu_d$ permutes without fixed points the
images of the generators of $F_n$. As a consequence, the intersection
of $\ker B$ with $\pi_1(X_{dr}^*,\mu_{dr})^{\mu_d}$ is trivial, and
$b$ is injective. The surjectivity of $b$ results from the surjectivity of $B$
and from the $\mu_d$-equivariance of the diagram.
\end{itemize}
\end{proof}

Together, the last two propositions prove the part of theorem
\ref{letheoreme} about symmetric groups, where
the only regular numbers are the divisors of $n$ and $n-1$.

Now consider the irreducible reflection group $G(p,1,n)$, 
denote by $X_{p,n}$ the complement of its discriminant.
As $G(p,1,n)$ is the centralizer of a $p$ regular
element of $\S_{pn}$, we can identify  $X_{p,n}$ and
$X_{pn}^{\mu_p}$ (in a way compatible with the action of $\mu_p$).

Let $d$ be a regular number for $G(p,1,n)$.
Consider the following diagram of inclusions:
$$\xymatrix{
X_{pn}^{\mu_e} \ar@{^{(}->}[rr] \ar@{^{(}->}[dr] & &
              X_{pn}^{\mu_p} \ar@{^{(}->}[dl] \\
& X_{pn} }$$
where $e$ is the lcm of $d$ and $p$. By proposition \ref{case1},
the $\pi_1$-images of downward arrows factorize through isomorphisms
with centralizers in $\pi_1(X_{pn},\mu_{pn})$ (all base points
being $\mu_{pn}$):
$$\xymatrix{
& \pi_1(X_{pn}^{\mu_e}) \ar@{^{(}->}[rr] \ar@{^{(}->}[dr] 
    \ar[dl]_\sim  & &
     \pi_1(X_{pn}^{\mu_p}) \ar@{^{(}->}[dl] \ar[dr]^\sim \\
\pi_1(X_{pn})^{\mu_e} \ar@{^{(}->}[rr] & & \pi_1(X_{pn})
& & \pi_1(X_{pn})^{\mu_p} \ar@{^{(}->}[ll] }$$      

Using the identification between $X_{p,n}$ and $X_{pn}^{\mu_p}$, the
part of the theorem about $G(p,1,n)$ says that the natural morphism
$\pi_1(X_{pn}^{\mu_e}) \rightarrow \pi_1(X_{pn}^{\mu_p}) $ induces an
isomorphism $\pi_1(X_{pn}^{\mu_e}) \stackrel{\sim}{\rightarrow}
 \pi_1(X_{pn}^{\mu_p})^{\mu_d}$; this is obvious on the above diagram,
which identifies $\pi_1(X_{pn}^{\mu_e})$ with $\pi_1(X_{pn})^{\mu_e}$
and $\pi_1(X_{pn}^{\mu_p})$ with $\pi_1(X_{pn})^{\mu_p}$.

To complete the proof of the theorem, it is enough to notice
that the mentioned exceptional groups are Shephard groups (see \cite{orso})
and that their discriminants are isomorphic to ones we have already
studied.

\end{document}